\documentclass{elsart}
\usepackage{epsfig}

\begin{document}

\begin{frontmatter}

\title{A conservative semi-Lagrangian method for oscillation-free computation 
of advection processes}
\author{Masato Ida}
\address{Satellite Venture Business Laboratory, Gunma University, \\
1--5--1 Tenjin-cho, Kiryu-shi, Gunma 376--8515, Japan\\
E-mail: ida@vbl.gunma-u.ac.jp}

\begin{abstract}
The semi-Lagrangian method using the hybrid-cubic-rational interpolation 
function [M. Ida, Comput. Fluid Dyn. J. 10 (2001) 159] is modified to a 
conservative method by incorporating the concept discussed in 
[R. Tanaka \textit{et al.}, Comput. Phys. Commun. 126 (2000) 232]. 
In the method due to Tanaka \textit{et al.}, 
not only a physical quantity but also its integrated quantity within a 
computational cell are used as dependent variables, and the mass conservation 
is completely achieved by giving a constraint to a forth-order polynomial used 
as an interpolation function. In the present method, a hybrid-cubic-rational 
function whose optimal mixing ratio was determined theoretically is employed 
for the interpolation, and its derivative is used for updating the physical 
quantity. The numerical oscillation appearing in results by the method due to 
Tanaka \textit{et al.} is sufficiently eliminated by the use of the hybrid 
function.
\end{abstract}

\begin{keyword}
Numerical method \sep Conservative method \sep Semi-Lagrangian \sep 
Interpolation \sep Cubic function \sep Rational function \sep 
Convexity preserving
\PACS 02.60.Cb \sep 02.70.-c \sep 47.11.+j
\end{keyword}

\end{frontmatter}

\section{Introduction}
The semi-Lagrangian method is a category of numerical methods for partial 
differential equations including an advection (or convection) term, and has 
been used mainly for atmospheric and geographic problems \cite{ref1,ref2}. In the semi-Lagrangian approach, the advection of a physical quantity is solved as 
a transportation and interpolation problem on fixed grids. The trajectory of 
the Lagrangian invariant that will reach a grid point is computed backwardly 
by the time integration of the advection velocity, and the invariant at the 
departure point, which is determined by an interpolation within a computational 
cell, is transported to the grid point. One can use a CFL-free time step in 
this approach. Since there are some well-established methods for calculating 
the trajectory such as the Runge-Kutta methods, recent studies on this 
approach are focused on the construction of interpolation functions.

One of problems of this approach is the difficulty of constructing a completely 
conservative scheme. There are only a few studies attacking to this problem. In 
Ref.~\cite{ref3} Priestley proposed a quasi-conservative semi-Lagrangian method 
by coupling a flux-corrected transport method and a high-order interpolation. 
In Ref.~\cite{add6}, Leonard \textit{et al.} proposed a conservative, explicit 
algorithm achieving stable high-CFL computations, by using an integral variable 
of a Lagrangian invariant as the dependent variable, and showed results for 
one-dimensional advection at constant velocity by coupling the algorithm with 
many interpolation schemes. In Ref.~\cite{add7}, Manson and Wallis modified the 
QUICKEST scheme \cite{add8} so that the stable high-CFL computation is 
achieved. Applications of their method to pure advection problems in 
non-uniform flow fields are shown in Ref.~\cite{add9}.
In Ref.~\cite{ref4} Lin and Rood discussed the use of readymade Eulerian 
schemes (the MUSCL \cite{add1} and PPM \cite{add2} schemes) modified 
to achieve stable computations even with a large time step.

Recently, Tanaka \textit{et al.} proposed a way to solve this difficulty. In 
Ref.~\cite{ref5,ref6}, they modified the CIP method, which is a 
non-conservative semi-Lagrangian method using the Hermite cubic interpolation 
function \cite{ref7,ref8}, into a conservative one by employing the 
cell-integrated quantity of a Lagrangian invariant as an additional dependent 
variable. A 4th- or 2nd-order polynomial is 
used as an interpolation function, and is constructed using a physical 
quantity, its first spatial derivative (only for the 4th-order one), and its 
cell-integrated quantity. The conservation of total mass is achieved by 
updating the cell-integral quantity in a conservation sense. The CIP method 
and its conservative variant using the 2nd-order polynomial (the CIP-CSL2 
method) are briefly reviewed in the next section.

The aim of this paper is to propose a variant of the CIP-CSL2 method by 
employing the hybrid-cubic-rational (HCR) interpolation function introduced 
recently by the author \cite{ref9}. The hybrid function is constructed with a 
combination of the cubic polynomial and a rational function, and is proposed 
to obtain oscillation-free results having higher resolution than those 
given by the rational method \cite{ref10} being a variant of the CIP. In the 
results using the CIP-CSL2 method one can sometimes observe numerical 
oscillation caused by the use of the classical polynomial \cite{ref6}. 
As is well known, the numerical oscillation becomes a serious 
problem for a certain examples such as the propagation of shock waves and the 
multiphase flow problems including large density jump. In Ref.~\cite{add5} 
the use of the rational function employed in the rational method is suggested 
to overcome the numerical oscillation. The author expects that combining 
the hybrid function with the concept of the CIP-CSL can derive an 
oscillation-free, conservative semi-Lagrangian method that has higher 
resolution than that by the rational method. In Sec.~\ref{secpre} we show that 
the hybrid method can easily be modified into a conservative one, and in 
Sec.~\ref{secexp} we demonstrate the accuracy and stability of the present 
method. In this paper, we discuss only the case of ${\rm CFL\;number} \le 1$; 
however, we believe that the present method is applicable to a high-CFL 
condition by incorporating the concept discussed in Ref.~\cite{ref6}.

\section{CIP method and its conservative variants}
The CIP method \cite{ref7,ref8} is a semi-Lagrangian solver for an advection 
equation,
\begin{equation}
\label{eq1}
{\frac{{\partial f}}{{\partial t}}} + u{\frac{{\partial f}}{{\partial x}}} = 
0{\rm ,}
\end{equation}
where $f$ is a dependent variable and $u$ is the advection velocity. The 
fact that the solution of this equation can be represented as
\begin{equation}
\label{eq2}
f(x,t + \Delta t) = f(X(x,t),t)
\end{equation}
allows us to solve Eq.~(\ref{eq1}) as an interpolation problem. Where $X$ is 
the trajectory of fluid particle that locates at $x$ at the time $t + \Delta t$,
\begin{equation}
\label{eq3}
X(x,t) = x + \int_t^{t - \Delta t} {u(X(x,\tau ),\tau )d\tau } .
\end{equation}
In the CIP method, the dependent variable $f$ is interpolated by the Hermite 
cubic function,
\[
C(X) = f_i^{}  + d_i^{} (X - x_i^{} ) + C2_i^{} (X - x_i^{} )^2  + C3_i^{} (X - x_i^{} )^3 ,
\]
constructed using both $f_i$ and its derivative $d_i$ (
$ = \partial f/\partial x|_i$) defined at each grid points $x_i$, where 
$i = 1,2, \cdots ,N$, $N$ is the number of grids, and $C2_i$ and $C3_i$ are 
determined from the condition of continuity \cite{ref7}. Equation (\ref{eq1}) 
and its spatial derivative, i.e.,
\begin{equation}
\label{eq4}
\frac{{\partial f(x,t + \Delta t)}}{{\partial x}} = \frac{{\partial X(x,t)}}{{\partial x}}\frac{{\partial f(X(x,t),t)}}{{\partial X}} ,
\end{equation}
are used to update $f_{i}$ and $d_{i}$, respectively. 
In the case where we assume the first-order accuracy in time for example, 
Eqs.~(\ref{eq2}) and (\ref{eq4}), respectively, are reduced to
\begin{equation}
\label{eq5}
f(x,t + \Delta t) = f(x - u(x,t)\Delta t,t) ,
\end{equation}
\begin{equation}
\label{eq6}
\frac{{\partial f(x,t + \Delta t)}}{{\partial x}} = \left[ {1 - \frac{{\partial u(x,t)}}{{\partial x}}\Delta t} \right]\frac{{\partial f(x - u(x,t)\Delta t,t)}}{{\partial X}} ,
\end{equation}
where $\Delta t$ is the time interval during a computational step. In the 
case of $u_i  < 0$, those equations are represented using the Hermite cubic 
function as
\begin{equation}
\label{eq7}
f_i^{n + 1}  = f_i^n  + d_i^n \xi  + C2_i^n \xi ^2  + C3_i^n \xi ^3 ,
\end{equation}
\begin{equation}
\label{eq8}
d_i^{n + 1}  = d_i^n  + 2C2_i^n \xi  + 3C3_i^n \xi ^2 ,
\end{equation}
where
\begin{equation}
\label{eq9}
C2_i^n  =  - (d_{i + 1}^n  + 2d_i^n  - 3S_{i + 1/2}^n )/h ,
\end{equation}
\begin{equation}
\label{eq10}
C3_i^n  = (d_{i + 1}^n  + d_i^n  - 2S_{i + 1/2}^n )/h^2 ,
\end{equation}
\[
S_{i + 1/2}^n  = (f_{i + 1}^n  - f_i^n )/h ,
\]
\[
\xi  =  - u_i \Delta t ,
\]
$h$ is the grid width, and the superscript $n$ denotes the time $t = n\Delta t$.

Basically, the CIP-CSL2 method is a solver for the equation of continuity,
\begin{equation}
\label{eq11}
\frac{{\partial f}}{{\partial t}} + \frac{{\partial uf}}{{\partial x}} = 0 .
\end{equation}
In this method, a quadratic function of
\begin{equation}
\label{eq12}
Q_i (X) = f_i  + 2\,q{\rm 1}_i (X - x_i ) + 3\,q2_i (X - x_i )^2
\end{equation}
is used to interpolate the physical quantity $f$, where the coefficients 
$q1_{i}$ and $q2_{i}$ are determined, for example in the case of $u_i  < 0$, 
from the constrains of
\[
Q_i (x_{i + 1} ) = f_{i + 1} \quad {\rm and} \quad \int_{x_i }^{x_{i + 1} } {Q_i (x)\,dx} /h = \rho _{i + 1/2},
\]
as
\begin{equation}
\label{eq13}
q1_i  =  - (f_{i + 1}  + 2f_i  - 3\rho _{i + 1/2} )/h ,
\end{equation}
\begin{equation}
\label{eq14}
q2_i  = (f_{i + 1}  + f_i  - 2\rho _{i + 1/2} )/h^2 ,
\end{equation}
where $\rho $ is the cell-integrated average of $f$ within a computational 
cell, defined at the cell centers, and is 
used as a dependent variable in this method (We should note that the 
definition of $\rho $ is slightly different from that in Ref.~\cite{ref5}). 
$\rho $ is updated so that the total mass is conserved. The 
mass flux flowing out from the cell [$x_i$, $x_{i + 1}$] through the point 
$x_i$ during $\Delta t$ can be calculated by
\[
\Delta \rho _i h = \int_{x_i }^{x_i  - u_i \Delta t} {Q_i (x)\,dx}  = f_i \xi  + q{\rm 1}_i \xi ^2  + q2_i \xi ^3 ,
\]
where we assume again the first-order accuracy in time. Using this quantity, 
$\rho $ is updated as follows:
\begin{equation}
\label{eq15}
\rho _{i + 1/2}^{n + 1}  = \rho _{i + 1/2}^n  + \Delta \rho _{i + 1}^n  - \Delta \rho _i^n .
\end{equation}
The use of this formula allows us the complete conservation of the total 
mass, i.e., the total sum of $\rho $.

The quantity $f$ is updated by a time splitting technique \cite{ref5}. We 
solve Eq.~(\ref{eq11}) by splitting into two phases as
\[
\frac{{\partial f}}{{\partial t}} + u\frac{{\partial f}}{{\partial x}} = 0 ,
\]
\[
\frac{{\partial f}}{{\partial t}} =  - \frac{{\partial u}}{{\partial x}}f .
\]
The former one shows the advection of $f$ and, therefore, is solved as an 
interpolation problem by using Eq.~(\ref{eq12}) as
\begin{equation}
\label{eq16}
f_i^{n + 1}  = f_i^n  + 2q{\rm 1}_i^n \xi  + 3q2_i^n \xi ^2 .
\end{equation}
The latter one represents the change of $f$ due to the compression or expansion 
and is solved in general by a finite difference technique \cite{ref5,ref7}.

As was pointed out already \cite{ref6}, the coefficients of the first- and 
second-order terms in Eq.~(\ref{eq16}) correspond to those in Eq.~(\ref{eq8}) 
by replacing $f$ and $\rho $ in Eq.~(\ref{eq16}) with $d$ and $S$, 
respectively. Furthermore, Eq.~(\ref{eq15}) also can be rewritten into a form 
like the conventional CIP. By introducing a variable $D$ defined as
\[
D_i  = \int_{x_0 }^{x_i } {f(x)\,dx} \;\;\left( { = \sum\limits_{j = 0}^{i - 1} {\int_{x_j }^{x_{j + 1} } {Q_j (x)\,dx} } } \right) ,
\]
one obtains
\begin{equation}
\label{eq17}
D_i^{n + 1}  = D_i^n  + \Delta \rho _i^n h = D_i^n  + f_i^n \xi  + q{\rm 1}_i^n \xi ^2  + q2_i^n \xi ^3 ,
\end{equation}
and
\[
\rho _{i + 1/2}^n  = (D_{i + 1}^n  - D_i^n )/h ,
\]
where we assume $\Delta \rho _0^n = 0$, namely no transfer of mass occur at 
the boundary $x_0$. Those equations reproduce Eq.~(\ref{eq15}), 
and replacements of $D$ and $f$ in Eq.~(\ref{eq17}) 
with $f$ and $d$, respectively, yield Eq.~(\ref{eq7}) (Note that the former 
replacement reduces $\rho $ to $S$). This result suggests us that the variants 
of the CIP method such as the rational method \cite{ref10} and the HCR method 
\cite{ref9} can be modified to a conservative one only by those replacements.

\section{A conservative method by the hybrid-cubic-rational interpolation}
\label{secpre}
In Ref.~\cite{ref9}, the author proposed a numerical solver for the advection 
equation by employing both the cubic \cite{ref7} and rational \cite{ref10} 
functions. In the method, the following combination of those functions is used 
as an interpolation function:
\begin{equation}
\label{eq18}
F(X) = \alpha \,R(X) + (1 - \alpha )\,C(X) ,
\end{equation}
where $R$ and $C$ show the rational and the cubic functions, respectively, 
and $\alpha$ denotes the weighting parameter whose range is $\alpha \in [0,1]$. 
Equation (\ref{eq18}) is reduced to the rational function for $\alpha = 1$ 
and to the cubic one for $\alpha = 0$. The value of $\alpha$ is determined 
theoretically so as to be the minimum of the values to which the 
convexity-preserving condition \cite{ref10} is satisfied. Where the 
convexity-preserving condition is expressed as 
that $F_{xx} (X) > 0$ for a concave data of $d_i  < S_{i + 1/2}  < d_{i + 1}$ 
or $F_{xx} (X) < 0$ for a convex data of $d_i  > S_{i + 1/2}  > d_{i + 1} $ is 
satisfied for the region of $X \in [x_i ,x_{i + 1} ]$. The resulting formulae 
are
\begin{equation}
\label{eq19}
F(k) = f_i  + d_i hk + (G1_i  + G2_i )k^2 ,
\end{equation}
\begin{equation}
\label{eq20}
\frac{{\partial F(k)}}{{\partial x}} = d_i  + [G1_i \frac{{Q_i  + D_i }}{{D_i }} + 2G2_i  + (1 - \alpha )(Q_i  - D_i )]\frac{k}{h} ,
\end{equation}
where
\[
G1_i  = \alpha P_i ^2 /D_i ,
\quad
G2_i  = (1 - \alpha )(2P_i  - D_i ) ,
\]
\[
D_i  = Q_i  + (P_i  - Q_i )k ,
\]
\[
\alpha  = \frac{{M_i (M_i  - 2)}}{{M_i (M_i  - 2) + 1}} ,
\]
\[
M_i  = \max [2,\;\max (\frac{{Q_i }}{{P_i }},\;\frac{{P_i }}{{Q_i }})]
\]
with
\[
P_i  = (S_{i + 1/2}  - d_i )h ,
\quad
Q_i  = (d_{i + 1}  - S_{i + 1/2} )h ,
\]
and
\[
k = \xi /h =  - u_i \Delta t/h .
\]
For the case of $u_i > 0$, we need replacements of $i + 1 \to i - 1$ and 
$h \to - h$. This interpolation function provides oscillation-free results, 
unlike the cubic one, and higher resolution of solution than that by the 
conventional rational function \cite{ref10}. Same as the conventional CIP, this 
method is constructed using a physical quantity $f_i$ and its first spatial 
derivative $d_i$; thus this method may be modified into a conservative method 
only by the replacements as discussed in the last section. Though the explicit 
definition and use of the cell-integrated average $\rho $ may be more 
convenient in an actual application, we now use $D$, instead of $\rho$, and 
adapt the replacements of
\begin{equation}
\label{eq21}
\left\{ {\begin{array}{*{20}c}
   {f \to D,}  \\
   {d \to f,}  \\
\end{array}} \right.
\end{equation}
for demonstration. Namely, the subroutine ``${\rm HCR}(f,d)$'' of the hybrid 
method whose parameters are $f$ and $d$ is used as a function of $D$ and $f$, 
i.e., ${\rm HCR}(D,f)$. The initial 
condition is set as follows. At first, $f$ is initialized using an arbitrary 
function $G(x)$, which expresses the initial profile of $f$, as
\[
f_i^0 = G(x_i ) .
\]
Next, the initial value of $D$ is calculated using a recurrence procedure of
\[
D_0^0 = 0 ,
\]
\begin{equation}
\label{eq22}
D_i^0  = D_{i - 1}^0  + h(f_i^0  + f_{i - 1}^0 )/2 ,
\quad
{\rm for}\;\;i = 1,2, \cdots ,N .
\end{equation}
Here the last term of Eq.~(\ref{eq22}) means the approximated integration of 
$f$ over a cell between $x_{i - 1}$ and $x_i$.

In the next section we show some numerical results using those procedures.

\section{Numerical experiments}
\label{secexp}
\subsection*{Example 1:}
\begin{figure}
\begin{center}
\leavevmode
\epsfxsize = 10 cm
\epsffile{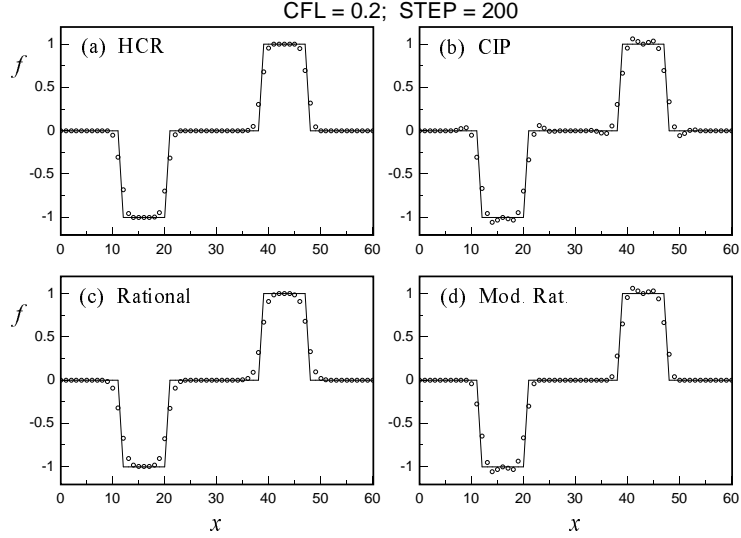}
\caption{
Linear propagation of square waves by (a): the HCR, (b): CIP, (c): rational, 
and (d): modified rational methods with the replacements of Eq.~(\ref{eq21}). 
Those results are with ${\rm CFL} = 0.2$ at $n = 200$. The circles and the 
solid lines show the numerical and theoretical results, respectively.}
\label{fig1}
\end{center}
\end{figure}

\begin{figure}
\begin{center}
\leavevmode
\epsfxsize = 10 cm
\epsffile{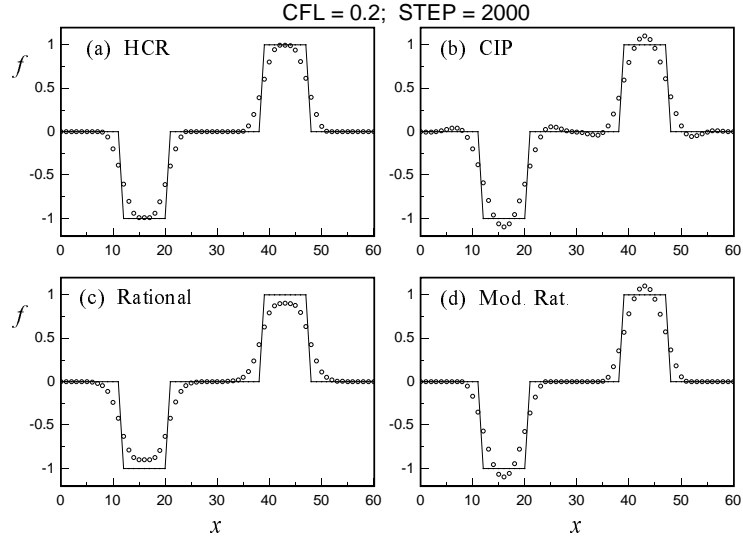}
\caption{
Same as Fig.~\ref{fig1} but $n = 2000$.}
\label{fig2}
\end{center}
\end{figure}

The first example is the linear propagation of square waves. In this example 
we use $u = 1$, $h = 1$ and
\[
G(x_i) = \left\{ {\begin{array}{*{20}c}
   { - 1,} \hfill & {{\rm for}\;\;13 \le i \le 21,} \hfill  \\
   1, \hfill & {{\rm for}\;\;40 \le i \le 48,} \hfill  \\
   0, \hfill & {{\rm elsewhere}{\rm .}} \hfill  \\
\end{array}} \right.
\]
Thus, the width of the square waves is $9h$. Figures \ref{fig1} and \ref{fig2} 
show the results using ${\rm CFL} = 0.2$ at $n = 200$ and $n = 2000$, 
respectively. Here we show results using four types of method, the HCR, CIP 
(the case of $\alpha = 0$), rational ($\alpha = 1$), and modified rational 
methods. In the modified rational method \cite{ref11}, a switching technique, 
which is represented as

\[
\gamma  = \left\{ {\begin{array}{*{20}c}
   {1,} \hfill & {d_i  \cdot d_{i + 1}  < 0,} \hfill  \\
   {0,} \hfill & {\rm otherwise,} \hfill  \\
\end{array}} \right.
\]
is added to the conventional rational method to select the interpolation 
function. The rational function is adopted in the case of $\gamma = 1$ and the 
cubic one in the case of $\gamma = 0$. It was pointed out 
that this technique sometimes breaks the preservation of the convexity of 
solution while this improves the dissipation property of the rational method 
\cite{ref9}. In the present example, this technique is adapted by using 
$f$, not by $d$, because of the replacements. In the figures, we see clearly 
that the HCR method provides most accurate results, which are less diffusive 
than those by the conventional rational method and less oscillatory than those 
by the CIP and 
modified rational methods. In Table 1, we show the calculated values at 
$n = 200$ around the point $i = 12$ at which the left 
discontinuity of the negative pulse locates. The calculated data by the CIP 
and modified rational methods are positive near the left side of the 
discontinuity and become minimums at $i = 14$, which are smaller than $-1$. 
On the contrary, the data by the HCR and conventional rational methods 
decrease monotonously as $i$ increases, and are not less than $-1$. 
Those results mean that only the HCR and conventional rational methods 
can keep the convexity of solution.

\begin{figure}
\begin{center}
\leavevmode
\epsfxsize = 10 cm
\epsffile{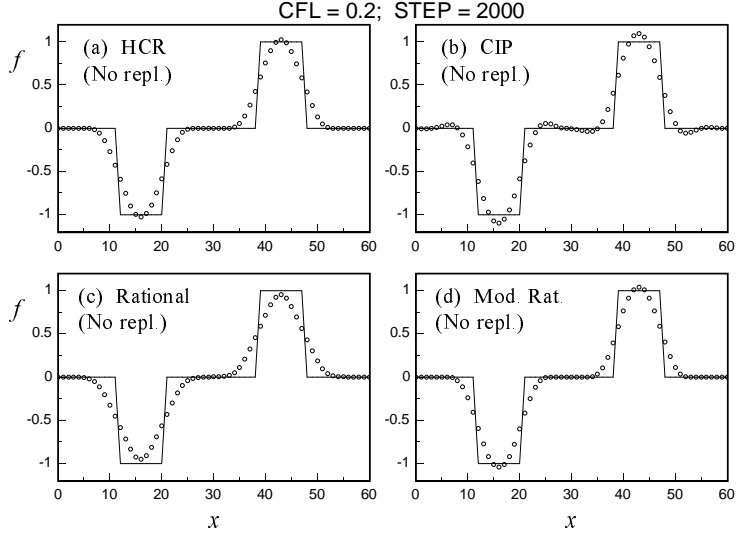}
\caption{
Same as Fig.~\ref{fig2} but the replacements of Eq.~(\ref{eq21}) are not 
adapted; namely those schemes are solved as non-conservative methods.}
\label{fig3}
\end{center}
\end{figure}

\begin{table}
\begin{center}
\begin{tabular}{|p{18pt}|p{70pt}|p{70pt}|p{70pt}|p{70pt}|}
\hline
$i$& 
Present& 
CIP& 
Rational& 
Mod. Rat. \\
\hline
4 \par 5 \par 6 \par 7 \par 8 \par 9 \par 10 \par 11 \par 12 \par 13 \par 14 \par 15 \par 16& 
0 \par 0 \par 0 \par -0.000001 \par -0.000044 \par -0.001716 \par -0.052075 \par -0.305191 \par -0.681895 \par -0.954887 \par -0.999656 \par -0.999996 \par -0.999997& 
-0.000014 \par -0.000986 \par -0.001887 \par 0.004413 \par 0.024674 \par 0.032729 \par -0.052964 \par -0.304522 \par -0.665011 \par -0.955764 \par -1.058063 \par -1.029841 \par -0.999959& 
0 \par -0.000004 \par -0.000037 \par -0.000345 \par -0.002738 \par -0.018069 \par -0.09232 \par -0.32068 \par -0.67146 \par -0.90682 \par -0.982846 \par -0.997188 \par -0.998634& 
0 \par 0 \par 0.000019 \par 0.000223 \par 0.000387 \par -0.000109 \par -0.04052 \par -0.277075 \par -0.647653 \par -0.951461 \par -1.059713 \par -1.031324 \par -1.000665 \\
\hline
\end{tabular}
\caption{Calculated values of the results corresponding to those shown in 
Fig.~\ref{fig1}.}
\end{center}
\end{table}

For comparison, we shall show results without the replacement of 
Eq.~(\ref{eq21}). Figure \ref{fig3} shows the results at $n = 2000$ given by 
using the same initial condition as the previous one. In this case, we set 
the initial value of the first derivative as $d_i^0  = 0$ at all points. As 
was pointed out in Ref.~\cite{ref6}, the CIP-CSL2 and CIP methods (the CIP 
with and without the replacements, respectively) give quite similar results 
each other. However, another methods provide obviously different results by 
whether the replacements are employed or not. It is interesting to point out 
that the result by the HCR method with the replacements is less diffusive than 
that by without the replacements, while we cannot explain this result for the 
moment.
\subsection*{Example 2:}
\begin{figure}
\begin{center}
\leavevmode
\epsfxsize = 10 cm
\epsffile{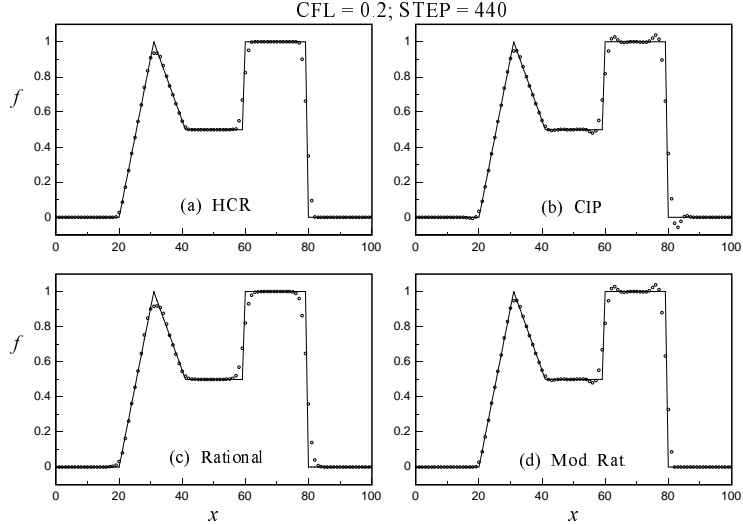}
\caption{
Linear propagation of triangular and square waves. Those results are with 
${\rm CFL} = 0.2$ at $n = 440$. The circles and the solid lines show the 
numerical and theoretical results, respectively.}
\label{fig4}
\end{center}
\end{figure}

\begin{figure}
\begin{center}
\leavevmode
\epsfxsize = 10 cm
\epsffile{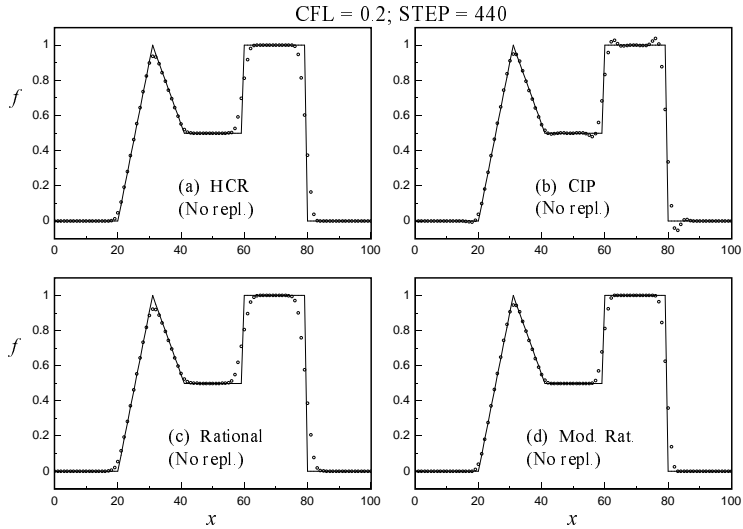}
\caption{
Same as Fig.~\ref{fig4} but the replacements are not adapted.}
\label{fig5}
\end{center}
\end{figure}

Next, we solve the problem used in Ref.~\cite{ref10} to compare the performance 
of the conventional rational method with those of the MUSCL \cite{add1} and PPM 
\cite{add2} schemes (See Fig.~1 in the paper). In this example, we use $h = 1$ 
and ${\rm CFL} = 0.2$. The initial value of $f$ is set to
\[
G(x_i = i) = \left\{ {\begin{array}{*{20}c}
   {(i - 20)/11,} \hfill & {{\rm for}\;\;20 \le i < 31,} \hfill  \\
   {1 - (i - 31)/20,} \hfill & {{\rm for}\;\;31 \le i < 41,} \hfill  \\
   {1/2,} \hfill & {{\rm for}\;\;41 \le i < 60,} \hfill  \\
   1, \hfill & {{\rm for}\;\;60 \le i < 80,} \hfill  \\
   0, \hfill & {{\rm elsewhere,}} \hfill  \\
\end{array}} \right.
\]
and the velocity, assumed as constant in time and space, is set to 
$u = 1$. Figures \ref{fig4} and \ref{fig5} show the results at $n = 440$, given 
by the four types of methods, with and without the 
replacements, respectively. In Ref.~\cite{ref10} it was concluded that the 
conventional rational method gives better representation of the sharp corner 
at the top of the triangular wave than those by the MUSCL and PPM methods. 
From the present result, however, we know that the HCR method and its 
conservative version achieve still better representation of that than the 
conventional rational method. The calculated maximum values of $f$ at the 
corner are 0.935 (conservative HCR), 0.937 (HCR), 0.916 (conservative 
rational), and 0.923 (rational). Furthermore, we know that the conservative HCR 
method gives almost equivalent (or slightly better) resolution of the square 
wave to that by the PPM, which gives the best one among those by the schemes 
discussed in Ref.~\cite{ref10}.

\subsection*{Example 3:}
Now we show results to a nonlinear problem in which the inviscid Burgers' 
equation,
\begin{equation}
\label{eq23}
\frac{{\partial u}}{{\partial t}} + u\frac{{\partial u}}{{\partial x}} = 0 ,
\end{equation}
is selected as the governing equation \cite{ref5,ref6}. This equation can be 
solved directly by the conservative methods discussed in this paper. For the 
calculation of $u$, we can use the procedures for $f$ because Eq.~(\ref{eq23}) 
is an advection equation (here $f$ means that after the replacements). 
For the calculation of $D$ ($= \int {u\,dx}$), we need to rewrite 
Eq.~(\ref{eq23}) into a conservation form as
\begin{figure}
\begin{center}
\leavevmode
\epsfxsize = 10 cm
\epsffile{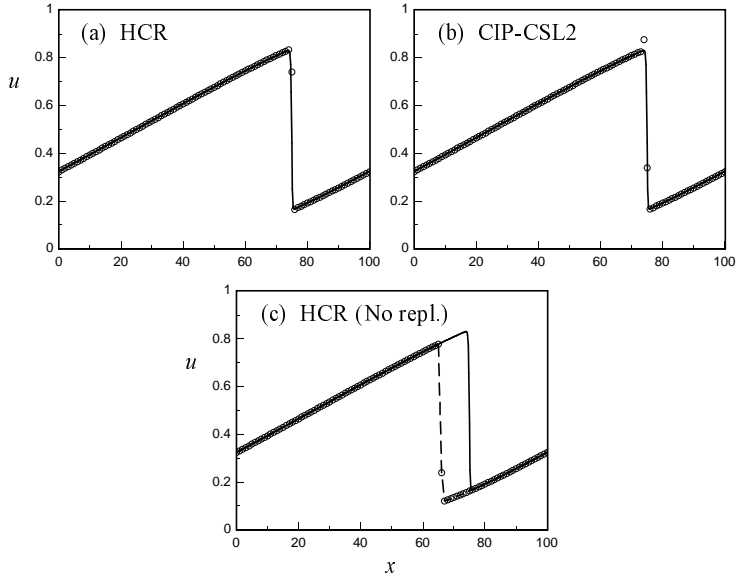}
\caption{
Result of inviscid Burgers' equation at $t = 100$ by (a): the HCR method with 
the replacements, (b): the CIP-CSL2 method, and (c): the HCR method without the 
replacements. The parameters are $h = 1$ and $\Delta t = 0.1$. The solid lines 
show results by the first-order upwind method with $h = 1/5$ and 
$\Delta t = 0.1/5$.}
\label{fig6}
\end{center}
\end{figure}

\begin{figure}
\begin{center}
\leavevmode
\epsfxsize = 10 cm
\epsffile{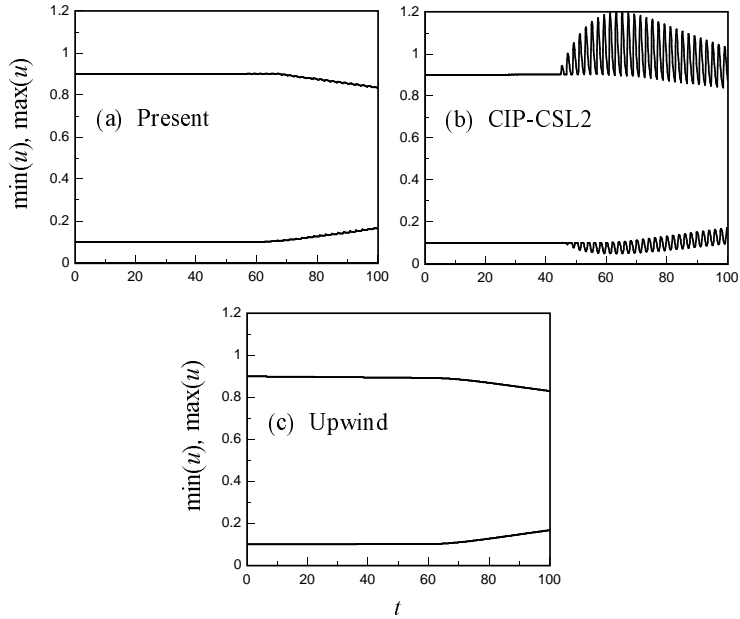}
\caption{
Maximum and minimum values in the results of inviscid Burgers' equation, as a 
function of time. Those results are by (a): the present and (b): CIP-CSL2 
methods with $h = 1$ and $\Delta t = 0.1$, and by (c): the first-order upwind 
method with $h = 1/5$ and $\Delta t = 0.1/5$.}
\label{fig7}
\end{center}
\end{figure}
\begin{equation}
\label{eq24}
\frac{{\partial u}}{{\partial t}} + \frac{{\partial (u^2 /2)}}{{\partial x}} = 0 .
\end{equation}
From this formula we found that the tranportation velocity $u$ used in the 
procedures for $D$ should be replaced with $u/2$ \cite{ref5}.

The initial condition and the grid width are set to
\[
G(x) = 0.5 + 0.4\cos (2\pi x/100)
\]
and $h = 1$. Figure \ref{fig6} shows results at $t = 100$ with 
$\Delta t = 0.1$. In this figure, the solid line shows the result by the 
first-order upwind method using the conservation formula (\ref{eq24}) with 
$h = 1/5$ and $\Delta t = 0.1/5$. Same as the result by the CIP-CSL2, that by 
the conservative HCR method represents the correct position of the shock wave. 
Figure \ref{fig6}(c) shows the result by the HCR method without the 
replacements. This result is obviously incorrect. Those results prove that the 
HCR method has gained complete conservation by the replacements. Figure 
\ref{fig7} shows the maximum and minimum values of $u$  as a function of time. 
In the result by the CIP-CSL2 we can observe unphysical strong oscillation due 
to the numerical dispersion of a method using a classical polynomial 
interpolation; however, such the oscillation cannot be seen in the result by 
the present method. This result allows us to consider that the present method 
is applicable to a shock problem even without any artificial viscosity, and is 
oscillation free even for the nonlinear problem.
\subsection*{Example 4:}
\begin{figure}
\begin{center}
\leavevmode
\epsfxsize = 10 cm
\epsffile{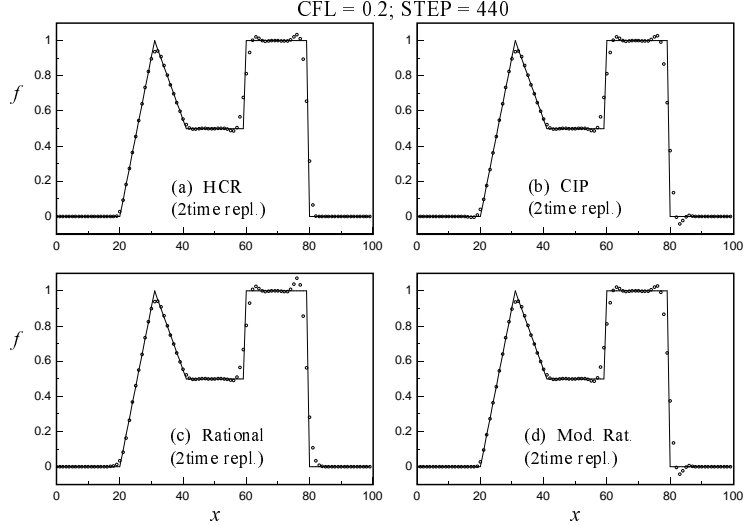}
\caption{
Same as Figs.~\ref{fig4} and \ref{fig5} but the two-time replacements are 
adapted.}
\label{fig8}
\end{center}
\end{figure}

Now we shall try an experimental study which uses further replacements. We 
already know that the CIP method provides an almost equivalent result 
whether the replacements (\ref{eq21}) are adapted or not. If the variables 
$D$ and $f$ are replaced with those integrated once more, what result is 
obtained? Figure \ref{fig8} shows results to \textit{Example 2} obtained by 
replacing $D_i^0$ and $f_i^0$ with
\[
E_i^0  \equiv \int_{x_0 }^{x_i } {D^0 \,dx}
\quad {\rm and} \quad
D_i^0
\]
and by updating those variables with the HCR and conventional methods. The 
parameters are same as those used in \textit{Example 2}. $f_i$ shown in 
Fig.~\ref{fig8} is given by $f_i = (E_{i + 1} - 2E_i + E_{i - 1})/h^2$. Despite 
the adaptation of two times of the replacements, the CIP method provides an 
almost equivalent result to those given without or with one-time 
replacements. However, one can observe some differences in the results using 
another methods from the previous ones. The numerical oscillation around the 
discontinuities of the square wave appears even in the results obtained by 
using the HCR and the conventional rational methods. In the result by 
the modified rational method, the oscillation can be observed not only at 
the top of the square wave but also at the bottom of it.

\section{An additional comment}
It is known that when no artificial viscosity is adapted, a numerical method 
for an advection equation based on the non-conservative formula cannot 
give correct propagation speed of shock in the region where the direction of 
the advection velocity changes spatially. In Ref.~\cite{ref5} the use of 
high-order polynomial that covers two cells are discussed for solving this 
problem. Now we present an alternative method to overcome this difficulty. 

Here we use again the Burgers' equation (\ref{eq23}) as the governing equation. 
Introducing $U \equiv u + C$, where $C$ is a constant, reduces Eq.~(\ref{eq23}) 
to
\[
\frac{{\partial U}}{{\partial t}} + U\frac{{\partial U}}{{\partial x}} - C\frac{{\partial U}}{{\partial x}} = 0 .
\]
This is also the Burgers' equation, but has a linear advection term (the last 
term of left hand side). We can solve this equation by splitting into two 
stages. First, we solve $\partial U/\partial t + U\partial U/\partial x = 0$ 
and update $U$ during the time interval $\Delta t$. Next, we shift the position 
of $U$ by $ - C\Delta t$. If $C$ is determined so that $U_i > 0$ or $U_i < 0$ 
for all $i$, we may be able to calculate the position of shock correctly. In 
the case where we set $C$ so that the condition of 
$m \left| {C\Delta t} \right| = h$ (where $m$ is a positive integer) is 
satisfied, the shift step can be solved only by one substitution like 
$U_{i + 1} = U_i$  per $m$ computational steps. This procedure is convenient 
since we need no additional interpolation function or difference equation; however, further discussions should be required for the practical application.

\section{Conclusion}
In this paper we have derived a conservative variant of the 
hybrid-cubic-rational method for an advection equation. It was proved by some 
numerical experiments that the hybrid method is oscillation free even when it 
is modified to a conservative one. Interestingly, the conservative hybrid 
method gives less-diffusive representation of the discontinuities of the square 
wave than that by the non-conservative one. We should try to make this result 
clear theoretically.

We expect that combining the present method with the hybrid 
interpolation-extrapolation method, which realizes the discontinuous 
representation of the density interface \cite{ref12}, yields a powerful tool 
for a semi-Lagrangian computation of multiphase flows like those treated in 
Refs.~\cite{add3,add4}.

\end{document}